\documentclass{article}
\usepackage{amsfonts,latexsym}
\usepackage{amsmath, times}
\usepackage{amssymb}
\usepackage{amsthm}
\usepackage{color}

\newcommand{\R}{\mathbb R}

\newcommand{\del}{\partial}
\newcommand{\e}{\varepsilon}
\newtheorem{theorem}{Theorem}[section]
\newtheorem{lemma}{Lemma}

\newtheorem{corollary}{Corollary}[section]
\newtheorem{definition}{Definition}[section]
\newtheorem{remark}{Remark}[section]
\newlength{\defbaselineskip}
\newcommand{\setlinespacing}[2]%
{\setlength{\baselineskip}{#1 \defbaselineskip}}

\makeatletter

\@addtoreset{equation}{section} \makeatother 
\begin{document}
\begin{center}
 {\Large   {Ground state  solutions  for  weighted biharmonic problem involving non linear exponential growth
}}
\end{center}
\vspace{0.2cm}

\begin{center}
 Brahim Dridi $^{1}$ and  Rached Jaidane $^{2}$

 \
\noindent\footnotesize   $^{1}$
Umm Al-Qura University, Faculty of Applied Sciences, Department of mathematics, P.O. Box $14035$,
Holly Makkah $21955$, Saudi Arabia.\\
Address e-mail: iodridi@uqu.edu.sa\\
\noindent\footnotesize $^2$ Department of Mathematics, Faculty of Science of Tunis, University of Tunis El Manar, Tunisia.\\
 Address e-mail: rachedjaidane@gmail.com\\
\end{center}

\vspace{0.5cm}
\noindent {\bf Abstract.}
In this article, we study the  following problem
$$\Delta(w(x)\Delta u) = \ f(x,u) \quad\mbox{ in }\quad B, \quad u=\frac{\partial u}{\partial n}=0 \quad\mbox{ on } \quad\partial B,$$
 where $B$ is the unit ball of $\mathbb{R}^{4}$ and  $ w(x)$
a singular weight of logarithm type. The reaction source
$f(x,u)$ is a radial function with respect to $x$ and it is critical in view of exponential inequality of Adams' type.
The existence result  is proved by using the constrained minimization in Nehari set coupled with the quantitative deformation lemma and degree theory results.\\

\noindent {\footnotesize\emph{Keywords:} Weighted Sobolev space, biharmonic operator, Critical exponential growth.\\
\noindent {\bf $2010$ Mathematics Subject classification}: $35$J$20$, $49$J$45$, $35$K$57$, $35$J$60$.}

\section{Introduction and Main results}
In  this paper, we consider  the fourth order weighted elliptic equation:
\begin{equation}\label{eq:1.1}
   \displaystyle (P)~~~~\left\{
      \begin{array}{rclll}
\Delta(w(x)\Delta u)  &=&  f(x,u)& \mbox{in} & B \\
        u&=\frac{\partial u}{\partial n}&0 &\mbox{on }&  \partial B,
      \end{array}
    \right.
\end{equation}
where $B=B(0,1)$ is the unit open ball in $\R^{4}$, $f(x,t)$ is a radial function with respect to $x$, the weight  $w(x)$ is given by \begin{equation}\label{eq:1.2}w(x)=\bigg(\log \frac{e}{|x|}\bigg)^{\beta},~~\beta\in(0,1)\cdot\end{equation}

In order to motivate our study, we first give a brief overview of the notion of exponential critical growth for problems of order superior or equal to $2$ in dimension $N\geq2$. We limit ourselves to Sobolev spaces $W^{1,N}_{0}(\Omega)$ and $W^{2,2}_{0}(\Omega)$ .We give some examples of applications as we go along.\\ In dimension $N\geq 2$, the critical exponential growth is given by the well known\
Trudinger-Moser inequality \cite{JMo,NST} $$\displaystyle\sup_{\int_{\Omega} |\nabla u|^{N}\leq1}\int_{\Omega}e^{\alpha|u|^{\frac{N}{N-1}}}dx<+\infty~~\mbox{if and only if}~~\alpha\leq \alpha_{N},$$
where $\alpha_{N}=\omega_{N-1}^{\frac{1}{N-1}}$  with $\omega_{N-1}$ is the area of the unit sphere $S^{N-1}$ in $\mathbb{R}^{N}$. This last result opened the way to study second order problems under nonlinearities with exponential growths and in non-weighted Sobolev spaces. For instance, we cite the following problem  in dimension $N=2$ \begin{equation*}-\Delta u=f(x,u)~~\mbox{in}~~\Omega\subset \mathbb{R}^{2}\end{equation*} which  have been studied considerably \cite{Adi1,FMR,LL1,MS}.\\

Later, the Trudinger-Moser inequality was improved to weighted inequalities \cite{ CR1, CR2}.
When the weight is of logarithmic type, Calanchi and Ruf \cite{CR3} extend the Trudinger-Moser inequality and proved the following results in the space $W_{0,rad}^{1,N}(B,\rho)=cl\{u \in
C_{0,rad}^{\infty}(B)~~|~~\int_{B}|\nabla u|^{N}\rho(x)dx <\infty\}$:
\begin{theorem}\cite{CR2} \label{th1.1}\begin{itemize}\item[$(i)$] ~~Let $\beta\in[0,1)$ and let $\rho$ given by $ \rho(x)=\big(\log \frac{1}{|x|}\big)^{\beta}$, then
  \begin{equation*}
 \int_{B} e^{|u|^{\gamma}} dx <+\infty, ~~\forall~~u\in W_{0,rad}^{1,N}(B,\rho),~~
  \mbox{if and only if}~~\gamma\leq \gamma_{N,\beta}=\frac{N}{(N-1)(1-\beta)}=\frac{N'}{1-\beta}
 \end{equation*}
and
 \begin{equation*}
 \sup_{\substack{u\in W_{0,rad}^{1,N}(B,\rho) \\ \int_{B}|\nabla u|^{N}w(x)dx\leq 1}}
 \int_{B}~e^{\alpha|u|^{\gamma_{N,\beta} }}dx < +\infty~~~~\Leftrightarrow~~~~ \alpha\leq \alpha_{N,\beta}=N[\omega^{\frac{1}{N-1}}_{N-1}(1-\beta)]^{\frac{1}{1-\beta}}
 \end{equation*}
where $\omega_{N-1}$ is the area of the unit sphere $S^{N-1}$ in $\R^{N}$ and $N'$ is the H$\ddot{o}$lder conjugate of $N$.
\item [$(ii)$] Let $\rho$ given by $\rho(x)=\big(\log \frac{e}{|x|}\big)^{N-1}$, then
  \begin{equation*}\label{eq:71.5}
 \int_{B}exp\{e^{|u|^{\frac{N}{N-1}}}\}dx <+\infty, ~~~~\forall~~u\in W_{0,rad}^{1,N}(B,\rho)
 \end{equation*} and
 \begin{equation*}\label{eq:71.6}
\sup_{\substack{u\in W_{0,rad}^{1,N}(B,\rho) \\  \|u\|_{\rho}\leq 1}}
 \int_{B}exp\{\beta e^{\omega_{N-1}^{\frac{1}{N-1}}|u|^{\frac{N}{N-1}}}\}dx < +\infty~~~~\Leftrightarrow~~~~ \beta\leq N.
 \end{equation*}
\end{itemize}
\end{theorem}
 These results paved the way for the study of second order weighted elliptic problems in dimension $N\geq2$. We point out that  recently, in the case, $V=0$ or $V\neq 0$, Baraket et al. \cite{BJ}, S. Deng, T. Hu and C-L. Tang and Calanchi et al. \cite{CRS,DHT}, have proved the existence of a nontrivial solution for
the following boundary value problem
$$
 \displaystyle \left\{
      \begin{array}{rclll}
    -\textmd{div} (\sigma(x)|\nabla u(x)|^{N-2}\nabla u(x) ) +V(x)|u|^{N-2}u &=& \ f(x,u)& \mbox{in} & B \\
        u&=&0 &\mbox{on }&  \partial B,
      \end{array}
         \right.
      $$
where $B$ is the unit ball in $\R^N,\; N\geq2,$ the radial positive weight $w(x)$ is of
logarithmic type, the function $f(x, u)$ is continuous in $B\times\R$ and behaves like $\exp\{e^{\alpha{t^{\frac{N}{N-1}}}}\}~~\mbox{as}~~t\rightarrow+\infty$, for some $\alpha>0$. The authors proved that there is a non-trivial solution to this problem using minimax techniques combined with Trudinger-Moser inequality.\\
   \\ 
   Now we will give an overview of  fourth order problems in relation to Adams' inegalities.\\
For bounded domains $\Omega\subset \mathbb{R}^{4}$, in \cite{Ada,BRFS} the authors extended the  Trudinger Moser inequality to the higher order space $ W_{0}^{2,2}(\Omega)$ and obtained the so called Adams' inequalities,
 \begin{equation*}
 \sup_{\substack{u\in S }}
 \int_{\Omega}~(e^{\alpha u^{2})}-1)dx < +\infty~~~~\Leftrightarrow~~~~ \alpha\leq 32 \pi^{2}
 \end{equation*}
 where $$S=\{u\in W^{2,2}_{0}(\Omega)~~|~~\displaystyle\big(\int_{\Omega}|\Delta u|^{2}dx \big)^{\frac{1}{2}}\leq 1\}.$$
  These results allowed to investigate fourth-order problems with subcritical or critical nonlinearity involving continuous potential  (see \cite{FS} , \cite{Chen}).\\
\begin{remark}The biharmonic equation in dimension $N>4$ \begin{equation*}\Delta^{2} u=f(x,u)~~\mbox{in}~~\Omega\subset \mathbb{R}^{N},\end{equation*}where the nonlinearity $f$ has subcritical and critical polynomial growth of power less than $\frac{N+4}{N-4}$, have been extensively studied \cite{BGT,EFJ,GGS} . \end{remark}
Before stating our results, let's start by defining our functional space.
 Let $\Omega \subset \R^{4}$, be a bounded domain and $w\in L^{1}(\Omega)$ be a nonnegative function. We introduce the weighted Sobolev space $$ W_{0}^{2,2}(\Omega,w)=cl\{u\in
C_{0}^{\infty}(\Omega)~~|~~\displaystyle\int_{B}w(x)|\Delta u|^{2}dx <\infty\}.$$
We will focus on radial functions on the unit ball $B$ and consider the subspace
$$\mathbf{E}= W_{0,rad}^{2,2}(B,w)=cl\{u\in
C_{0,rad}^{\infty}(B)~~|~~\displaystyle\int_{B}w(x)|\Delta u|^{2}dx <\infty\},$$
 endowed with the norm $$\|u\|=\displaystyle\big(\int_{B}w(x)|\Delta u|^{2}dx\big)^{\frac{1}{2}}.$$
  We note that this norm is issued from the product scalar
$$\langle u,v\rangle=\int_{B} w(x)~~ \Delta u \Delta v~dx.$$

  The choice of the weight induced in (\ref{eq:1.2}) and the space $\mathbf{E}$ are also motivated by the following exponential inequality.
\begin{theorem}\cite{WZ} \label{th1.2} ~~Let $\beta\in(0,1)$ and let $w$ given by (\ref{eq:1.2}), then

 \begin{equation}\label{eq:1.3}
 \sup_{\substack{u\in W_{0,rad}^{2,2}(B,w) \\  \int_{B}|\Delta u|^{2}w(x)dx \leq 1}}
 \int_{B}~e^{\displaystyle\alpha|u|^{\frac{2}{1-\beta} }}dx < +\infty~~~~\Leftrightarrow~~~~ \alpha\leq \alpha_{\beta}=4[8\pi^{2}(1-\beta)]^{\frac{1}{1-\beta}}
 \end{equation}

 \end{theorem}

Let $\gamma:=\displaystyle\frac{2}{1-\beta}$,  in view  of inequality (\ref{eq:1.3}), we say that $f$ has critical growth at $+\infty$ if there exists some $\alpha_{0}>0$,
\begin{equation}\label{eq:1.4}
\lim_{s\rightarrow +\infty}\frac{|f(x,s)|}{e^{\alpha s^{\gamma}}}=0,~~~\forall~\alpha ~~\mbox{such that}~~ \alpha_{0}<\alpha\leq \alpha_{\beta}~~
\mbox{and}~~~~\lim_{s\rightarrow +\infty}\frac{|f(x,s)|}{e^{\alpha s^{\gamma}}}=+\infty,~~\forall~\alpha<\alpha_{0}\leq \alpha_{\beta}.\end{equation}
Now, let's state our assumptions. In this paper, we deal  with problem $(P)$ under critical growth nonlinearities.
Furthermore, we suppose that $f(x,t)$ satisfies the following hypothesis:
\begin{enumerate}
\item[$(V_{1})$] $f: B \times \mathbb{R}\rightarrow\mathbb{R}$ is continuous  and radial in $x$.
\item[$(V_{2})$] There exist $\theta > 2$  such that  we have
 $$0 < \theta F(x,t)\leq tf(x,t), \forall (x, t)\in~~B\times\mathbb{R} \setminus\{0\} $$
  where
$$F(x,t)=\displaystyle\int_{0}^{t}f(x,s)ds.$$
\item [$(V_{3})$]  For each $x\in B$,~$\displaystyle t \mapsto \frac{f(x,t)}{|t|}~~\mbox{is increasing for}~~ t\in~~\mathbb{R} \setminus\{0\}$.
 \item [$(V_{4})$]   $\displaystyle\lim_{t\rightarrow 0}\frac{|f(x,t)|}{|t|}=0.$
 \item[$(V_{5})$] There exist $p$ such that $2<\theta <p $ and $C_p > 1$ such that
$$ sgn(t) f(x,t) \geq C_p \vert t\vert^{p-1}, \quad \mbox{for all} \; (x,t)\in B \times \R,$$
where $sgn(t) = 1$ if $t > 0,$ $sgn(t) = 0$ if $t = 0,$ and $sgn(t) = -1$ if $t < 0.$
\end{enumerate}
We give an example of such nonlinearity. The nonlinearity $f(x,t)=C_{p}|t|^{p-2}t+|t|^{p-2}t\exp(\alpha_{0} |t|^{\gamma})$ satisfies the assumptions $(V_{1})$, $(V_{2})$, $(V_{3})$ ,$(V_{4})$ and $(V_{5})$.\\


We will consider the following definition of  solutions.
\begin{definition}\label {def1.1}
We say that a function $u\in \mathbf{E}$ is a weak solution of the problem $(P)$ if
\begin{equation*}\label {eq:1.9}
\int_{B}w(x)~\Delta u .\Delta \varphi~  dx =
\int_{B}f(x,u) \varphi dx,~~\forall~\varphi \in \mathbf{E}.
\end{equation*}
\end{definition}
Let $\mathcal{J} :\mathbf{E} \rightarrow \R$ be the functional given by
 \begin{equation}\label{eq:1.5}
\mathcal{J}(u)=\frac{1}{2}\int_{B}w(x)~|\Delta u|^{2}dx-\int_{B}F(x,u)dx,
\end{equation}
where
$$ F(x,t)=\displaystyle\int_{0}^{t}f(x,s)ds.$$
The energy functional $\mathcal{J}$ is well defined and of class $C^{1}$ since there exist $ a,~ C>0$ positive constants and there exists $t_{1} >1$ such that
\begin{equation*}\label {eq:1.11}
|f(x,t)|\leq C e^{a ~t^{\gamma}}, ~~~~~~\forall |t| >t_{1},
\end{equation*}
whenever the nonlinearity $f(x,t)$ is critical at $+\infty$.\\It is standard to check that critical points of $\mathcal{J}$ are precisely weak solutions of $(P)$. Moreover, we have $$\langle\mathcal{J}'(u),\varphi\rangle=\mathcal{J'}(u)\varphi=\int_{B}\big(w(x)~ \Delta u~\Delta \varphi\big)~dx-\int_{B}f(x,u)~ \varphi~dx~,~~\varphi \in\mathbf{E}\cdot$$Our strategy consists in finding  solutions which minimize the corresponding energy
functional $\mathcal{J} $ among the set of all solutions to problem $(P)$. To this end, we define the
 Nehari set as:

$$\displaystyle\ \mathcal{N}:=\{u\in \mathbf{E}:\langle \mathcal{J}'(u),u\rangle=0, u\neq 0\}.$$

In other words, we try to find a minimize of the energy functional $\mathcal{J}$ over the following minimization problem,
\begin{equation*}\label {eq:1.12}
\displaystyle m=\inf_{u\in\mathcal{N}}\mathcal{J}(u)
\end{equation*}

To our best knowledge, there are no results for  solutions to the weighted biharmonic  equation with critical exponential nonlinearity on the weighted Sobolev space $\mathbf{E}$.\\
Now, we give our main result as follows:
\begin{theorem}\label{th1.3}~~ Let $f(x,t)$ be a function that has a critical growth at
$+\infty$. Suppose that   $(V_{1})$, $(V_{2})$, $(V_{3})$, $(V_{4})$ and $(V_{5})$  are satisfied.
 The  problem $(P)$ has a radial solution with minimal energy provided
\begin{equation}\label{eq:1.6} C_{p} >\max\bigm\{1, \bigg(\frac{2\theta~~ m_{p}}{(\theta-2)}\big(\frac{2(\alpha_{0}+\delta)}{\alpha_{\beta}}\big)^{1-\beta}\bigg)^{\frac{p-2}{2}}\bigm\}
\end{equation}
 where $\delta >0$, $m_p = \displaystyle\inf_{u\in\mathcal{N}_p} J_p(u)>0$,
 $$J_p(u) :=\frac{1}{2}\|u\|^{2} -\frac{1}{p}\int_{B} \vert u\vert^p dx  $$ and
 $$\mathcal{N}_p:= \{ u \in \mathbf{E}, u \neq 0\; \mbox{and}\; \langle J_p'(u),u\rangle =0\}.$$
 \end{theorem}

In general the study of  fourth order partial differential equations is considered an interesting topic. The interest in studying such equations was stimulated by their applications in micro-electro-mechanical systems, phase field models of multi-phase systems, thin film theory,surface diffusion on solids, interface dynamics, flow in Hele-Shaw cells, see \cite{D, FW, M}. However  many applications are generated
by the weighted elliptic problems, such as  the study of traveling waves in suspension bridges, radar imaging  (see, for example \cite{AEG, LL}).\\

This present work is organized as follows: in Section $2$, we present some necessary
preliminary knowledge about functional space and some preliminaries results. In section $3$, we give  some technical key lemmas . In section $4$, we study an auxiliary problem which will be of great use to prove our main result. Section $5$ is devoted to the proof of the Theorem \ref{th1.3}.\\
Finally,  we note that a constant $C$ may change from line to another and sometimes we index the constants in order to show how they change. Also, we shall use the notation $\vert u \vert_p$ for the norm in the Lesbegue space $L^p(B)$.

\section{Weighted Sobolev Space setting and preliminaries results }
Let $\Omega \subset \R^{N}$, $N\geq2$,  be a bounded domain in $\R^{N}$ and let $w\in L^{1}(\Omega)$ be a nonnegative function. To deal with weighted operator, we need to introduce some functional spaces $L^{p}(\Omega,w)$, $W^{m,p}(\Omega,w)$, $W_{0}^{m,p}(\Omega,w)$ and some of their properties that will be used later. Let $S(\Omega)$ be the set of all measurable real-valued functions defined on $\Omega$ and two measurable functions are considered as the same element if they are equal almost everywhere.\\\\
Following  Drabek et al. and Kufner in \cite{DKN}, the weighted Lebesgue space $L^{p}(\Omega,w)$ is defined as follows:
$$L^{p}(\Omega,w)=\{u:\Omega\rightarrow \R ~\mbox{measurable};~~\int_{\Omega}w(x)|u|^{p}~dx<\infty\}$$
for any real number $1\leq p<\infty$.\\
This is a normed vector space equipped with the norm
$$\|u\|_{p,w}=\Big(\int_{\Omega}w(x)|u|^{p}~dx\Big)^{\frac{1}{p}}.$$
For $w(x)=1$, one finds the standard  Lebesgue space $L^{p}(\Omega)$ endowed with the norm $|u|_{p}=\Big(\int_{\Omega}|u|^{p}~dx\Big)^{\frac{1}{p}}.$\\
For $m\geq 2$, let $w$ be a given family of weight functions $w_{\tau}, ~~|\tau|\leq m,$ $$w=\{w_{\tau}(x)~~x\in\Omega,~~|\tau|\leq m\}.$$
In \cite{DKN}, the  corresponding weighted Sobolev space was  defined as
$$ W^{m,p}(\Omega,w)=\{ u \in L^{p}(\Omega)~ ~\mbox{such that}~~ D^{\tau} u \in L^{p}(\Omega)~\mbox{for all}~~|\tau|\leq m-1,~~D^{\tau} u \in L^{p}(\Omega,w)~\mbox{for all}~~|\tau|= m\} $$
endowed with the following norm:

\begin{equation*}\label{eq:2.2}
\|u\|_{W^{m,p}(\Omega,w)}=\bigg(\sum_{ |\tau|\leq m-1}\int_{\Omega}|D^{\tau}u|^{p}dx+\displaystyle \sum_{ |\tau|= m}\int_{\Omega}|D^{\tau}u|^{p}w(x) dx\bigg)^{\frac{1}{p}}.
\end{equation*}

If we suppose also that $w(x)\in L^{1}_{loc}(\Omega)$, then $C^{\infty}_{0}(\Omega)$ is a subset of $W^{m,p}(\Omega,w)$ and we can introduce the space $$W^{m,p}_{0}(\Omega,w)$$
as the closure of $C^{\infty}_{0}(\Omega)$ in $W^{m,p}(\Omega,w).$ Moreover, the following embedding is compact  $$W^{m,p}(\Omega,w)\hookrightarrow\hookrightarrow W^{m-1,p}(\Omega)\cdot$$
Also, $(L^{p}(\Omega,w),\|\cdot\|_{p,w})$ and $(W^{m,p}(\Omega,w),\|\cdot\|_{W^{m,p}(\Omega,w)})$ are separable, reflexive Banach spaces provided that $w(x)^{\frac{-1}{p-1}} \in L^{1}_{loc}(\Omega)$.\\
 Then the space $\mathbf{E}$ is a Banach and reflexive. The space $\mathbf{E}$ is endowed with the norm
$$ \| u\|=\left(\int_{B}w(x)|\Delta u|^{{2}}dx \right )^{\frac{1}{2}}$$
which is equivalent to the following norm (see lemma \ref{lemr} bellow) $$\|u\|_{W_{0,rad}^{2,2}(B,w)}=\displaystyle\big(\int_{B}u^{2}dx +\int_{B} |\nabla u|^{2}~dx+\int_{B}w(x)|\Delta u|^{2}dx \big)^{\frac{1}{2}}\cdot$$
We also have the continuous embedding   $$\mathbf{E}\hookrightarrow L^{q}(B)~~\mbox{for all}~~q\geq 1.$$
 Moreover, $\mathbf{E}$ is compactly
embedded in $L^{q}(B)$  for all $q \geq1$ . In fact, we have



\begin{lemma}\label{lemr}

\item[(i)] Let $u$ be a radially symmetric
 function in $C_{0}^{2}(B)$. Then, we have\begin{itemize}\item[$(i)$]\cite{WZ}
 $$|u(x)|\leq \displaystyle\frac{1}{2\sqrt{2}\pi}\frac{||\log(\frac{e}{|x|})|^{1-\beta}-1|^{\frac{1}{2}}}{\sqrt{1-\beta}}\displaystyle \int_B w(x)|\Delta u|^2dx \leq \displaystyle\frac{1}{2\sqrt{2}\pi}\frac{||\log(\frac{e}{|x|})|^{1-\beta}-1|^{\frac{1}{2}}}{\sqrt{1-\beta}}\|u\|^{2}\cdot$$

\item[(ii)] The norms $\|.\| $ and $\|u\|_{ W_{0,rad}^{2,2}(B,w)}=\displaystyle\big(\int_{B}u^{2}dx +\int_{B} |\nabla u|^{2}~dx+\int_{B}|\Delta u|^{2}w(x)dx \big)^{\frac{1}{2}}$ are equivalents.
\item[(iii)]  The following embedding
is continuous $$\mathbf{E}\hookrightarrow L^{q}(B)~~\mbox{for all}~~q\geq 1.$$
\item[(iv)]$\mathbf{E}$ is compactly
embedded in $L^{q}(B)$  for all $q \geq2$ .
\end{itemize}
\end{lemma}
\textit{Proof }

$(i)$ see \cite{WZ}\\

 $(ii)$~~ By
 Poincar\'{e} inequality,  for all $u \in W_{0,rad}^{1,2}(B)$
 $$\int_B  \displaystyle | u|^{2}
 \leq C \int_B  \displaystyle |\nabla u|^{2}.
 $$ Using the Green formula, we get
$$
\int_B  \displaystyle |\nabla u|^{2}= \int_B \nabla u \nabla u  = - \displaystyle \int_B u \Delta
 u + \underbrace{\displaystyle\int_{\del B} u \frac{\partial u}{\partial n}}_{= 0} ~
 \leq  \displaystyle \Big| \int_B u \Delta u \Big|\cdot
 $$
By  Young inequality, we get for all $\varepsilon>0$
 $$
 \displaystyle \Big| \int_B u \Delta u \Big| \leq
  \displaystyle \frac{1}{2 \varepsilon} \displaystyle \int_B |\Delta u|^2 + \displaystyle \frac{\varepsilon}{2} \displaystyle
 \int_B |u|^2\leq  \displaystyle \frac{1}{2 \varepsilon} \displaystyle \int_B w(x)|\Delta u|^2 + \displaystyle \frac{\varepsilon}{2} \displaystyle
 \int_B |u|^2.$$
 Hence
 $$ (1 - \displaystyle \frac{\varepsilon}{2} C^2) \displaystyle \int_B |\nabla u |^{2} ~ \leq ~
 \displaystyle \frac{1}{2\varepsilon} \displaystyle\int_ B  w(x)|\Delta u |^{2},$$
  then, $$\displaystyle\int_{B}u^{2}dx +\int_{B} |\nabla u|^{2}~dx+\int_{B} w(x)|\Delta u|^{2}dx \leq C \int_{B} w(x)|\Delta u|^{2}dx\leq C\|u\|^{2}\cdot$$
 Then $(ii)$ follows.\\
 $(iii)$ and $(iv)$. Since $w(x)\geq1$, then the following embedding are continuous and compact   $$\mathbf{E}\hookrightarrow W^{2,2}_{0,rad}(B,w)\hookrightarrow W^{2,2}_{0,rad}(B)\hookrightarrow L^{q}(B)~~\forall q\geq 2$$ and from $(i)$, we have that $\mathbf{E} \hookrightarrow L^{q}(B)$ is continuous for all $q\geq1$.
This concludes the lemma.\hfill $\Box$\\

\section{ Some technical lemmas }
We begin by some key lemmas.\\
In the following we assume, unless otherwise stated, that the function $f$ satisfies the conditions $(V_{1})$ to $(V_{4})$. Let $u\in\mathbf{E}$ with $u\not\equiv 0 $ a.e. in the ball $ B$, and we deﬁne the function
$\displaystyle\Upsilon_{u} : [0, \infty)  \rightarrow\R$
 as
\begin{equation}\label{eq:3.1}
 \displaystyle\Upsilon_{u}(t) = \mathcal{J}(tu ).
\end{equation}
It's clear that $\displaystyle\Upsilon'_{u}(t)=0$ is equivalent to $tu\in \mathcal{N}$.
\begin{lemma}\label{lem1}\begin{itemize}
\item [(i)]For  each  $u \in\mathbf{E}$ with  $u\neq0$  ,  there  exists  an  unique $t_{u}>0$,
 such that
$\displaystyle t_{u}u \in\mathcal{N}.$ In particular,  the set $\mathcal{N}$ is nonempty and $\mathcal{J}(u)>0$, for every $u\in \mathcal{N}$.
\item[(ii)]  For all $t \geq 0 $ with $t \neq t_{u},$ we have
$$\displaystyle \mathcal{J}(tu ) <  \mathcal{J}(t_{u}u)\cdot$$\end{itemize}

\end{lemma}
Proof.$(i)$ \\
Since $f$ is  critical, and from $(V_{1})$ and $(V_{4})$,  for all $\varepsilon> 0$, there exist  positive constants $C_{1} = C_{1}(\epsilon )$  and $C'_{1} = C'_{1}(\epsilon )$ such that
\begin{equation}\label{eq:3.2}
f(x,t)t\leq \varepsilon |t|^{2} +C_{1}|t|^{s }\exp(\alpha|t|^{\gamma}  )\mbox{ for all }	\alpha > \alpha_{0}, s > 2.
\end{equation}
and \begin{equation}\label{eq:3.3}
F(x,t)\leq  \frac{1}{2}\varepsilon |t|^{2} +C'_{1}|t|^{s }\exp(\alpha|t|^{\gamma}  )\mbox{ for all }	\alpha > \alpha_{0}, s > 2.
\end{equation}
 Now, given $u\in\mathbf{E}$ fixed with  $u \neq0$ . From (\ref{eq:3.3}), for all $\varepsilon>0$, we have
\begin{align}\label{eq:3.4}
 \displaystyle\Upsilon_{u}(t) = \mathcal{J}(tu )
 & =\frac{1}{2}t^{2}\|u\|^{2}-\int_{B}F(x,tu)tudx\nonumber   \\
 & \geq \frac{1}{2}t^{2}\|u\|^{2} \nonumber
 -\frac{1}{2}\epsilon t^{2} \int_{B}|u|^{2}dx\nonumber - C'_{1}\int_{B}|tu|^{s}\exp (\alpha t|u|^{\gamma})dx
 \end{align}
   Using the
 H\"{o}lder inequality, with $a, a' > 1$ such that $\displaystyle\frac{1}{a} + \frac{1 }{a'}  = 1$, and Sobolev embedding Lemma \ref{lemr}, we get
   \begin{align}
 \displaystyle\Upsilon_{u}(t) & \geq \frac{1}{2}t^{2}\|u\|^{2} - C_{2}\frac{1}{2}\epsilon t^{2}\|u\|^{2}\nonumber
   - C_{1}\left( \int_{B}|tu|^{a's}dx\right)^{ \frac{1 }{a'} }\left(\int_{B}\exp (\alpha t a|u|^{\gamma})dx\right)^{\frac{1}{a}}\nonumber \\ &\geq \left (\frac{1}{2}-\frac{1}{2}\epsilon C_{2} \right)\|tu\|^{2}- \left(\int_{B}\exp \big(\alpha  a\|tu\|^{\gamma}\big(\frac{|u|}{\|u\|}\big)^{\gamma}\big)dx\right)^{\frac{1}{a}} C_{3}\|tu\|^{s}\nonumber
   \end{align}
  By (\ref{eq:1.3}), the last integral is finte provided $t>0$ is chosen small enough such that  $\displaystyle\alpha a\|tu\|^{\gamma}\leq \alpha_{\beta}$. Then,
\begin{align*}
\displaystyle\Upsilon_{u}(t) & \geq \left(\frac{1}{2}-\frac{1}{2}\epsilon C_{2} \right)\|tu\|^{2}- C_{4}\|tu\|^{s}~~\mbox{with}~~\alpha a\|tu\|^{\gamma}\leq \alpha_{\beta}~~\mbox{and}~~\alpha >\alpha_{0}\end{align*}
holds. Choosing  $\epsilon > 0$ such that $\displaystyle\frac{1}{2}-\frac{1}{2}\epsilon C_{2} > 0$ and since $ s > 2$, we obtain,\begin{equation}\label{eq:3.4}  \displaystyle\Upsilon_{u}(t) > 0 ~~ \mbox{for small}~~ t>0 .\end{equation}

Now, From $(V_{3})$, we can derive that there exist $C_{5}, C_{6} > 0 $ such that
\begin{equation}\label{eq:3.5}
F (x,t)\geq  C_{5}|t|^{\theta}-C_{6}.
\end{equation}
Then, by using $(\ref{eq:3.5})$, we get
\begin{align*}
 \displaystyle\Upsilon_{u}(t)
  &= \mathcal{J}(tu ) \nonumber  \\
  & \leq \frac{1}{2} t^{2}\|u\|^{2}-C'_{5}|t|^{\theta}\|u\|^{\theta}
  - C_{6}|B|
\end{align*}
Since $\theta>2$, we get that \begin{align}\label{eq:3.6}\displaystyle\Upsilon_{u}(t)\rightarrow -\infty~~\mbox{as}~~t\rightarrow+\infty .\end{align} Hence, from (\ref{eq:3.4}) and (\ref{eq:3.6}),  there exists at least one $t_{u} > 0$ such that $\Upsilon'_{u}(t_{u})=0$, i.e. $t_{u} u\in \mathcal{N}$.  \\

Now we will show the uniqueness of  $t_{u}$. Let $s>0$ such that $su\in \mathcal{N}$. Then we get

 $\displaystyle\langle\mathcal{J}'(t_{u}u),  t_{u}u \rangle= 0,$
 $\displaystyle\langle\mathcal{J}'(su), su \rangle= 0,$ and

\begin{equation}\label{eq:3.7}
 \displaystyle \|su\|^{2} =\int_{B}f(x,su)su dx
\end{equation}

\begin{equation} \label{eq:3.8}
  \|t_{u}u\|^{2}=\int_{B}f(x,t_{u}u)tudx
\end{equation}

Combining (\ref{eq:3.7}) and (\ref{eq:3.8}), we deduce that
 $$ 0=\int_{B}\frac{f(x,t_{u}u)}{t_{u}u}u^{2}dx-\int_{B}\frac{f(x,su)}{su}u^{2}dx.$$
 It follows from $(V_{4})$  that $t\mapsto \frac{f(x,t)}{t}$ is increasing for $t > 0$, which implies
 that $t_{u} = s$.  This completes the proof of $(i)$.\\
$(ii)$ Follows from $(i)$ , since $\mathcal{J}(t_{u}u)=\displaystyle\max_{t\geq 0}\Upsilon_{u}(t)$. \\In the sequel, we prove that sequences in $\mathcal{N}$ cannot converge to $0$.
\begin{lemma}\label{lem2.2}
	Assume that  $(V_1)-(V_4)$ hold. Then for any $u \in \mathbf{E}$ with $u \neq 0$ such that
	$\langle \mathcal{J}'(u),u\rangle\leq0$, the unique maximum point  of  $\Upsilon_{u}$ on $\mathbb{R}_+ $ satisfies $0< t_u\leq 1$.
\end{lemma}

\noindent Proof:\\
	
	Since $t_u u\in\mathcal{N}$, we have
	\begin{equation} \label{eq:3.9}
	\begin{aligned}
t_u^2 \Vert u\Vert^2 = \int_{B}f(x,t_u u) t_u u dx.
	\end{aligned}
	\end{equation}
	Furthermore, since $\langle J'(u),u\rangle\leq0$, we have
	$$     \Vert u\Vert^2 \leq  \int_{B}f(x, u) u dx.$$
	Then by  \eqref{eq:3.9}, we have
	\begin{equation} \label{eq:3.10}
	\begin{aligned}
	(t_u^{-2}-1)\Vert u\Vert^2
	\geq \int_{B}\Big(\frac{f(x,t_u u)}{t_u u}-\frac{f(x, u)}{u}\Big)u^2\,dx.
	\end{aligned}
	\end{equation}
	Obviously, the left hand side of \eqref{eq:3.10} is negative for
	$t_u>1$ whereas the right hand side is positive, which is a
	contradiction. Therefore $0<t_u\leq 1$.

\begin{lemma}\label{lem3}
 For all $u\in \mathcal{N}$,
 \begin{itemize}
   \item [$(i)$] there exists $\kappa>0$ such that\\
   $ \|u\|  \geq \kappa ;$
   \item[$(ii)$] $\mathcal{J}(u) \geq (\frac{1}{2}-\frac{1}{\theta})\|u\|^{2}$
 \end{itemize}
 \end{lemma}
  Proof. $(i)$
We argue by contradiction.  Suppose that there exists a sequence $\{u_{n}\} \subset \mathcal{N} $
such that $u_{n}\rightarrow 0$ in $ \mathbf{E}.$
Since $\{u_{n}\} \subset \mathcal{N}$, then  $\displaystyle\langle\mathcal{J}'(u_{n}) ,u_{n} \rangle =0$.
Hence, it follows from (\ref{eq:3.2}), (\ref{eq:3.3}) and the radial Lemma \ref{lemr} that

\begin{align}\label{eq:3.11}
\|u_{n}\|^{2}&=\int_{B}f(x,u_{n})u_{n}dx \\
 & \leq \epsilon  \int_{B}|u_{n}|^{2}dx +  C_{1}\int_{B}|u_{n}|^{s}\exp(\alpha |u_{n}|^{\gamma} )dx\nonumber\\
 & \leq \epsilon  C_{6} \|u_{n}\|^{2} +  C_{1}\int_{B}|u_{n}|^{s}\exp(\alpha |u_{n}|^{\gamma} )dx\nonumber
\end{align}
Let $a>1$ with $\frac{1}{a}+\frac{1}{a'}=1$. Since $u_{n}\rightarrow 0\mbox{ in }~~\mathbf{E}$,
for $n$ large enough, we get
$\displaystyle\|u_{n}\|\leq(\frac{\alpha_{ \beta}}{\alpha a})^{\frac{1}{\gamma}}$. From H\"{o}lder inequality, (\ref{eq:1.3}) and again the radial Lemma \ref{lemr},
 we have
\begin{align*}
  \int_{B}|u_{n}|^{s}\exp(\alpha |u_{n}|^{\gamma} )dx&\leq
  \left( \int_{B}|u_{n}|^{sa'}dx\right)^{\frac{1}{a'}}
  \left(\int_{B}\exp \big(\alpha  a\|u^{+}\|^{\gamma}\big(\frac{|u^{+}|}{\|u^{+}\|}\big)^{\gamma}\big)dx\right)^{\frac{1}{a}} \\
  &\leq C_{7} \left( \int_{B}|u_{n}|^{sa'}dx\right)^{\frac{1}{a'}}\leq C_{8} \|u_{n}\|^{s}
\end{align*}
Combining (\ref{eq:3.11}) with the last inequality, for $n$ large enough, we obtain

\begin{equation}\label{eq:3.12}
 \|u_{n}\|^{2}\leq \epsilon C_{6} \|u_{n}\|^{2}
 + C_{8} \|u_{n}\|^{s}
\end{equation}
Choose suitable $\epsilon > 0 $ such that  $1-\epsilon C_{6} > 0$.
Since  $2 < s$, then (\ref{eq:3.12}) contradicts the fact that $u_{n}\rightarrow 0\mbox{ in }\mathbf{E}$.\\
$(ii)$ Given $u \in \mathcal{N}$, by the deﬁnition of $\mathcal{N}$ and $(V_{3})$ we obtain
\begin{align*}
  \mathcal{J}(u) &= \mathcal{J}(u)-\frac{1}{\theta}\langle\mathcal{J}'(u), u\rangle \\
  &= \frac{1}{2}\|u_{n}\|^{2}-\frac{1}{\theta}\|u_{n}\|^{2}
  + \big(\int_{B}\frac{1}{\theta}f(x,u)u-F(x,u)dx\big)\\
  &\geq (\frac{1}{2}-\frac{1}{\theta})\|u\|^{2}
\end{align*}
Lemma \ref{lem3} implies that $\mathcal{J}(u) > 0$  for all $u\in \mathcal{N}$.
As a consequence, $\mathcal{J}$ is bounded by below in $\mathcal{N}$, and therefore
$\displaystyle m:=\inf_{u\in \mathcal{N}} \mathcal{J}(u) $
 is well-deﬁned.\\ In the following lemma we prove that if the minimum of $\mathcal{J}$ on $\mathcal{N}$ is achieved in some $u\in\mathcal{N}$, then $u$  is a critical point of $\mathcal{J}$.

\begin{lemma}\label{lem5}
If $ u_{0}\in \mathcal{N} $ satisfies $\mathcal{J}(u_{0})=m$, then $\displaystyle\mathcal{J}'(u_{0})=0.$
\end{lemma}
\noindent Proof:We argue by contradiction. We assume that $\displaystyle\mathcal{J}'(u_{0})\neq 0$. By the continuity of $\mathcal{J}'_{\lambda}$,
there exist $\iota, \delta\geq 0$ such that
\begin{equation}\label{eq:3.13}
\displaystyle \|\mathcal{J}'_{\lambda}(v)\|_{\mathbf{E}^{\ast}}\geq\iota \mbox{ for all }~~v ~~\mbox{such that}~~ \|v-u_{0}\|\leq \delta.
\end{equation}
 Let $\displaystyle D=\left(1-\tau,1+\tau\right)\subset \mathbb{R}$ with $\displaystyle \tau \in(0, \frac{\delta}{4\|u_{0}\|})$ and define $g:D\rightarrow\mathbf{E}$, by

$$\displaystyle g(\rho)=\rho u_{0},  \rho \in D$$
By virtue of $u_{0} \in \mathcal{N}$, $\mathcal{J}(u_{0})=m$ and Lemma \ref{lem1}, it is easy to see that
\begin{equation}\label{eq:3.14}
\displaystyle\bar{m}:=\max_{\partial D} \mathcal{J}\circ g<m~~ \mbox{and}~~ \mathcal{J}(g(\rho))<m,~~\forall ~~\rho\neq 1.
\end{equation}
 Let $\epsilon:=\min\{\frac{m-\bar{m}}{2}, \frac{\iota\delta}{16}\}$, $S_{r}:=B(u_{0},r),r\geq0$
and $\displaystyle\mathcal{J}^{a}:=\mathcal{J}^{-1}(]-\infty,a]).$
 According to the quantitative deformation Lemma $[\cite{Wi}, \mbox{ Lemma }2.3]$,
there exists a deformation $\eta \in C\left(\mathbf{E}, \mathbf{E}\right)$ such that:
 \begin{itemize}
   \item [$(1)$] $\eta( v)=v,$ if $v\not\in \mathcal{J}^{-1}([ m-\epsilon,m+\epsilon])\cap S_{\delta}$
   \item[$(2)$] $\eta\left( \mathcal{J}^{m+\epsilon }\cap S_{\frac{\delta}{2}}\right)\subset \mathcal{J}^{m-\epsilon}$,
    \item[$(3)$] $\mathcal{J}(\eta( v))\leq \mathcal{J}(v)$, for all $v\in  \mathbf{E}. $
 \end{itemize}
By lemma \ref{lem1} $(ii)$, we have $\mathcal{J}( g(\rho))\leq m$. In addition, we have,$$\|g(\rho)-u_{0}\|=\|(\rho-1)u_{0}\|\leq \frac{\delta}{4},~~\forall \rho\in D\cdot$$ Then, $g(\rho)\in S_{\frac{\delta}{2}}$ for $\rho\in \bar{D}$.
Therefore, it follows from $(2)$ that
\begin{equation}\label{eq:3.15}
\max_{\rho\in \bar{D}}\mathcal{J}(\eta( g(\rho))\leq m-\epsilon.
\end{equation}
In the following, we prove that $\eta( g(D))\cap\mathcal{N}$ is nonempty. And in this case it contradicts (\ref{eq:3.15}) due to the definition
of $m$.
 To do this, we first define
$$\bar{g}(\rho):=\eta( g(\rho)),$$
\begin{align*}
  \Upsilon_{0}(\rho)& = \langle \mathcal{J}'(g(\rho)), u_{0} \rangle,
\end{align*}
and
$$\Upsilon_{1}(\rho):=(\frac{1}{\rho}\langle \mathcal{J}'(\bar{g}(\rho),(\bar{g}(\rho)) \rangle.$$
We have that for $\rho\in \overline{D}$, $\mathcal{J}(g(\rho))\leq \overline{m}<m-\varepsilon$. Then, $\bar{g}(\rho)=\eta (g(\rho))=\rho u_{0}$. Hence,
\begin{align}\label{eq:3.16}
  \Upsilon_{0}(\rho)& =  \Upsilon_{1}(\rho), \forall \rho \in \overline{D}
\end{align}
On one hand, we have that $\rho=1$ is the unique critical point of $\Upsilon_{0}$. So by degree theory, we get that $d^{0}(\Upsilon_{0},\mathcal{J},0)=1$. On the other hand, from (\ref{eq:3.16}), we deduce that  $d^{0}(\Upsilon_{1},\mathcal{J},0)=1$. Consequently, there exists $\overline{\rho}\in D$ such that $\overline{g}(\overline{\rho})\in \mathcal{N}$. This implies that
$$m\leq \mathcal{J}(\overline{g}(\overline{\rho}))=\mathcal{J}(\eta(g(\overline{\rho})).$$
This contradicts (\ref{eq:3.15}) and finish the proof of the Lemma.
\section{Auxilary problem}
In this section, in order to prove our existence  result , we consider
the auxiliary problem \begin{equation}\label{pau}
\displaystyle \left\{
\begin{array}{rclll}
 \Delta (w (x) \Delta u) &=&  \displaystyle \vert u\vert^{p-2}u& \mbox{in} & B \\
u&=\frac{\partial u}{\partial n}&=0 &\mbox{on }&  \partial B,
\end{array}
\right.
\end{equation}
where $p$ is the constant that appear in the hypothesis $(V_5)$. We have associated to problem \eqref{pau} the functional
$$J_p(u) :=\frac{1}{2}\|u\|^{2} -\frac{1}{p}\int_{B} \vert u\vert^p dx  $$ and the Nehari manifold
$$\mathcal{N}_p:= \{ u \in \mathbf{E}, u \neq 0\; \mbox{and}\; \langle J_p'(u),u \rangle = 0\}.$$
Let $m_p = \inf_{\mathcal{N}_p} J_p(u)>0$, we have the following results for $J_{p}$.
\begin{lemma}\label{ljp1} Given $u\in \mathbf{E}, u\neq 0$, there exists a unique $t>0$ such that $t u\in\mathcal{N}_{p}$. In addition, $t$ satisfies \begin{equation}\label {JP1}J_{p}(tu)=\displaystyle \max_{s\geq0}J_{p}(su)\end{equation}.
\end{lemma}
\noindent Proof: Let $\displaystyle\gamma(s) = J_{p}(su )=\frac{1}{2}s^{2}\|u\|^{2}-\frac{s^{p}}{p}|u|^{p}_{p}$, for $s>0$. Since $p>2$, we have that $\gamma(s)>0$ for $s>0$ small enough and $\gamma(s)\rightarrow -\infty~~\mbox{as}~~~s\rightarrow -\infty$. Hence, there exists a $t>0$ satisfying (\ref{JP1}). In particular, $tu\in \mathcal{N}_{p}$. Moreover, $\gamma'(t)=0$ if and only if $t=\big(\frac{\|u\|^{2}}{|u|^{p}_{p}}\big)^{\frac{1}{p-2}}$.\\
As a consequence, we have
\begin{corollary} Let $u\in \mathbf{E}, u\neq 0$. Then  $u\in \mathcal{N}_{p}$ if and only if $J_{p}(tu)=\displaystyle \max_{s\geq0}J_{p}(su)$.
\end{corollary}
\begin{lemma}\label{ljp2}
 For all $u\in \mathcal{N}_{p}$,
 \begin{itemize}
   \item [$(i)$] there exists $\kappa_{0}>0$ such that\\
   $ \|u\|  \geq \kappa_{0} ;$
   \item[$(ii)$] $\mathcal{J}_{p}(u) \geq (\frac{1}{2}-\frac{1}{p})|u|^{p}_{p}$
 \end{itemize}
 \end{lemma}
\begin{lemma} \label{ljp3}There exists $w_{p} \in \mathcal{N}_{p}$ such that $J_{p}(w_{p}) = m_{p}$ and
$m_{p}=\displaystyle\frac{p-2}{2p}|w_{p}|^{p}_{p}.$
\end{lemma}
\noindent Proof:
	Let sequence $(w_n) \subset \mathcal{N}_{p} $ satisfy $\displaystyle \lim_{n \rightarrow +\infty} J_{p}(w_n) = m_{p}$. It is clearly that $(w_n)$ is bounded by Lemma \ref{ljp2}. Then, up to a subsequence, there exists $w_{p} \in \mathbf{E}$ such that
	
	\begin{equation}\label{eq:4.3}\begin{array}{ll}
	w_{n} \rightharpoonup w_{p}~~~~&\mbox{in}~~\mathbf{E},\\
	w_{n} \rightarrow w_{p}~&\mbox{in}~~L^{q}(B),~~\forall q\geq 2,\\
	w_{n} \rightarrow w_{p} ~~&\mbox{a.e. in }~~B.
	\end{array}
	\end{equation}
	We claim that $w_{p}\ne 0$. Suppose, by contradiction, $w_{p}= 0$. From the definition
	of $\mathcal{N}_{p}$ and \eqref{eq:4.3} , we have that $ \displaystyle\lim_{n \rightarrow +\infty} \Vert w_n \Vert^2 =0$, which contradicts Lemma \ref{ljp2}. Hence, $w_{p}\ne 0$ .
	
	From the lower semi continuity of norm and \eqref{eq:4.3}, it follows that
 	\begin{equation}\label{eq:4.4}
 	\Vert w_{p}\Vert^2 \leq  \liminf_{n \rightarrow +\infty}(\Vert w_n\Vert^2 )
 	\end{equation}	

 On the other hand, by using $\langle J'(w_n),w_n\rangle =0$ and (\ref{eq:4.3}), we have

 	\begin{equation}\label{eq:4.5} \liminf_{n \rightarrow +\infty}(\Vert w_n\Vert^2 ) = \liminf_{n \rightarrow +\infty} \int_{B}| w_n|^{p} dx=\int_{B} |w_{p}|^{p} dx.\end{equation}
 From \eqref{eq:4.4} and \eqref{eq:4.5} we deduce that $\langle J'(w_{p}),w_{p}\rangle \leq 0$.
 Then, as in Lemma \ref{lem2.2} this implies that there exists $s_u \in (0, 1] $ such that $s_uw_{p}\in\mathcal{N}_{p}$. Thus, by the lower semi continuity of norm and (\ref{eq:4.3}), we get that

$$\begin{array}{rclll}
\displaystyle m_{p} \leq J_{p}(s_uw_{p})&=&J(s_uw_{p}) -\displaystyle\frac{1}{2} \langle J_{p}'(s_uw_{p}),s_uw_{p}\rangle\\&=&\displaystyle\big(\frac{1}{2}-\frac{1}{p}\big)s^{p}_{u}\int_{B}|w_{p}|^{p}dx
\\ &\leq& J_{p}(w_{p}) -\displaystyle\frac{1}{2} \langle J_{p}'(w_{p}),w_{p}\rangle\\
&=&\displaystyle\frac{1}{2}\|w_{p}\|^{2} -\frac{1}{p}\int_{B}|w_{p}|^{p}dx-\frac{1}{2}\|w_{p}\|^{2}+\frac{1}{2}\int_{B}|w_{p}|^{p}dx\\
&\leq&\displaystyle\liminf_{n\rightarrow +\infty}\Big[\displaystyle\frac{1}{2}\|w_{n}\|^{2} -\frac{1}{p}\int_{B}|w_{n}|^{p}dx-\frac{1}{2}\|w_{n}\|^{2}+\frac{1}{2}\int_{B}|w_{n}|^{p}dx\Big]\\
&\leq&\displaystyle \liminf_{n\rightarrow+\infty}\big[J_{p}(w_n) -\displaystyle\frac{1}{2} \langle J_{p}'(w_n), w_n\rangle\big] = m_{p}.
\end{array}$$	
Therefore, we get that $J_{p}(w_{p})=m_{p}$, which is the desired conclusion. \\
\section{Proof of Theorem 1.2}
Now, we will obtain an important estimate for the level $m$. That will be a powerful tool in order
to obtain an appropriate bound of the norm of a minimizing sequence for $m$ in $\mathcal{N}$.
	\begin{lemma}\label{lem9}
	Assume that $(V_1)-(V_5)$ and \eqref{eq:1.6} are satisfied. It holds that
		\begin{equation} \label{eq:5.1} m \leq \frac{\theta-2}{2\theta}\bigg(\frac{ \alpha_{\beta}}{2 (\alpha_0 + \delta)}\bigg)^{1-\beta}.\end{equation}
\end{lemma}

\noindent Proof:
From Lemma \ref{ljp3}, there exists $w_{p} \in \mathcal{N}_p$ such that $J_p(w_{p}) = m_p$ and $J_p'(w_{p})=0.$ Consequently, we get	
\begin{equation}\label{eq:5.2}
\frac{1}{2}\Vert w_{p}\Vert^2 - \frac{1}{p}\int_{B} \vert w_{p}\vert^p ~dx =  m_p
\end{equation}
and
	\begin{equation} \label{eq:5.3}
\Vert w_{p}\Vert^2 = \int_{B} \vert w_{p}\vert^p ~dx.
\end{equation}

By virtue of $(V_5)$ and \eqref{eq:5.3}, we have $\langle \mathcal{J}'(w_{p}),w_{p}\rangle \leq 0$ which together with Lemma \ref{lem2.2} yielding that there is a unique  $s \in (0, 1] $ such that $sw_{p}\in\mathcal{N}$. Using $(V_5)$, \eqref{eq:5.2} and  \eqref{eq:5.3}, we obtain
$$\begin{array}{rclll}
\displaystyle m &\leq& \mathcal{J}(sw_{p})
\\  &\leq& \displaystyle\frac{s^2}{2}\|w_{p}\|^2
\displaystyle - \frac{C_p s^p}{p}\vert w_{p}\vert_p^p \\  &=& \displaystyle\big(\frac{s^2}{2} -  \frac{C_p s^p}{p}\big)\vert w_{p}\vert_p^p  \\
 &\leq&\displaystyle \max_{\xi >0}(\frac{\xi^2}{2} - C_p\frac{\xi^p}{p}) |w|_p^p
\end{array}$$
By some straightforward algebraic manipulations, we get
\begin{equation}\label{eq:5.4}
m \leq  C_p^{\frac{-2}{p-2}}\frac{p-2}{2p}|w_{p}|_p^p.
\end{equation}
Note that by using \eqref{eq:5.2}, \eqref{eq:5.3} and the fact that $p>\theta>2$, we have
\begin{equation}\label{eq:5.5}
(\frac{1}{2}- \frac{1}{p})|w_{p}|_p^p  =\frac{1}{2}\|w_{p}\|^{2}- \frac{1}{p}|w_{p}|_p^p= m_p.
\end{equation}
Thus, by combining \eqref{eq:5.4} and \eqref{eq:5.5}, we obtain
\begin{equation}\label{eq:5.6}
m \leq  C_p^{\frac{-2}{p-2}} m_p.
\end{equation}
Therefore, by \eqref{eq:1.6} and \eqref{eq:5.6}, we obtain that \eqref{eq:5.1} holds.
The following result gives us some compactness properties of minimizing sequences.

	\begin{lemma}\label{lemma9}
	If $(u_n) \subset \mathcal{N} $ is a minimizing sequence for $m$, then there exists $u \in \mathbf{E}$ such that
	
	$$\int_{B}f(x,u_{n})u_n dx \rightarrow \int_{B}f(x,u)u dx $$
	and
	$$\int_{B}F(x,u_{n}) dx \rightarrow \int_{B}F(x,u)dx. $$
\end{lemma}

\noindent Proof:
	We must prove the first limit, since the second one is analogous.  It is sufficient to prove that $g(u_n(x))$ is convergent in $L^1(B)$, where $g(u_n(x))$ is defined by \begin{equation*}\label{eq4.1}
	\int_{B}f(x, u_n)~ u_n dx \leq \e
	\int_{B}| u_n|^2\,dx + C \int_{B}|u_n|^{q} exp(\alpha (u_n)^{\gamma})dx = g(u_n(x)), \quad \mbox{for all}\; \alpha > \alpha_0\, \mbox{and}\, q>2.
	\end{equation*}
First note that
	\begin{equation}\label{eq:5.7}
	|u_{n}|^q \rightarrow |u|^q~~ \mbox{in}~~ L^{2}(B).
	\end{equation}	
	On the other hand, by $(V_{2})$, we obtain that
		\begin{equation} \label{eq:5.8}
	\begin{aligned}
	\displaystyle m=\limsup_{n \rightarrow +\infty}  \mathcal{J}(u_n)
	& = \limsup_{n \rightarrow +\infty} \big(\mathcal{J}(u_n)-\frac{1}{\theta}\langle \mathcal{J}'(u_n),u_n\rangle \big) \\
	& =\limsup_{n \rightarrow +\infty} \big(  \frac{\theta -2}{2\theta}\| u_n\|^{2} +\frac{1}{\theta}\int_{B}\big(f(x, u_n) u_n- \theta F(x, u_n)\big)dx\big) \\
	& >  \frac{\theta -2}{2\theta}\limsup_{n \rightarrow +\infty} \|u_n\|^2,
	\end{aligned}
	\end{equation}
	which together with Lemma \ref{lem9} gives that $\displaystyle \limsup_{n \rightarrow +\infty} \|u_n\|^\gamma < \frac{\alpha_{\beta}}{2(\alpha_0+\delta)}\cdot$
	
    Now choosing $\alpha= \alpha_0+\delta$, we have from Theorem \ref{th1.2} that
	\begin{equation}\label{eq:5.9}
	\int_{B} exp(2\alpha |u_n|^{\gamma})dx \leq \int_{B}\exp\big(2(\alpha_{0} + \delta )\|u_{n}\|^{\gamma}(\frac{u_{n}}{\|u_{n}\|})^{\gamma}\big)dx\leq \int_{B}\exp\big(\alpha_{\beta}(\frac{u_{n}}{\|u_{n}\|})^{\gamma}\big)dx.
	\end{equation}
	Then it follows by (\ref{eq:1.3}) that there is $M > 0$ such that
		\begin{equation}\label{eq:5.10}
 \int_{B} exp(2 \alpha |u_n|^{\gamma})dx \leq M.
	\end{equation}
	Since
	\begin{equation}\label{eq:5.11}
	exp(\alpha |u_n|^{\gamma}) \rightarrow exp(\alpha |u|^{\gamma})\:  \mbox{a.e in }~~B.
	\end{equation}
	From \eqref{eq:5.10} and [\cite{Ka}, Lemma 4.8], we get that
	\begin{equation}\label{eq:5.12}
	exp(\alpha |u_n|^{\gamma}) \rightharpoonup exp(\alpha |u|^{\gamma})\:  \mbox{ in }~~L^{2}(B).
	\end{equation}
	Now using \eqref{eq:5.7}, \eqref{eq:5.12} and [\cite{Ka}, Lemma 4.8] again, we conclude that
	\begin{equation}\label{eq4.6}
	\int_{B} f(x, u_n~) u_n dx \rightarrow \int_{B} f(x, u)~ u dx.
	\end{equation}
In the sequel, we give an important result:
\begin{lemma}\label{lem11} Assume that the conditions ($V_{1}$), ($V_{2}$) and ($V_{3}$) are satisfied. Then, for each $x\in B$, we have
	$$tf(x, t) -2 F(x, t)~~ \mbox{is increasing for} ~ t> 0 ~ \mbox{ and decreasing for} ~ t< 0.  $$
	In particular, $tf(x, t) - 2 F(x, t)> 0~~\mbox{for all} ~~(x,t) \in B \times \R\setminus \{0\}.$
\end{lemma}
\noindent Proof: Assume that $0<t<s$. For each $x\in B$, we have
$$\begin{array}{rlll}
\displaystyle
tf(x,t)-2F(x,t)&=& \displaystyle\frac{f(x,t)}{t}t^{2}-2F(x,s)+2\int^{s}_{t}f(x,\nu)d\nu \\
&<&\displaystyle \frac{f(x,t)}{s}t^{2}-2F(x,s)+\frac{f(x,s)}{s}(s^{2}-t^{2})\\

&=&sf(x,s)-2 F(x,s)\cdot
\end{array}$$
The proof in the case $t<s<0$ is similar.\\
 The assertion $tf(x, t) - 2 F(x, t)> 0~~\mbox{for all} ~~(x,t) \in B \times \R\setminus \{0\}$ comes from ($V_{2}$).\hfill $\Box$\\

\begin{lemma} \label{lem12}There exists $w_{0} \in \mathcal{N}$ such that $\mathcal{J}(w_{0}) = m$ .
\end{lemma}
\noindent Proof:
	Let sequence $(w_n) \subset \mathcal{N} $ satisfy $\displaystyle \lim_{n \rightarrow +\infty}\mathcal{J} (w_n) = m$. It is clearly that $(w_n)$ is bounded by Lemma \ref{ljp2}. Then, up to a subsequence, there exists $w_{0} \in \mathbf{E}$ such that
	
	\begin{equation}\label{eq:5.14}\begin{array}{ll}
	w_{n} \rightharpoonup w_{0}~~~~&\mbox{in}~~\mathbf{E},\\
	w_{n} \rightarrow w_{0}~&\mbox{in}~~L^{q}(B),~~\forall q\geq 2,\\
	w_{n} \rightarrow w_{0} ~~&\mbox{a.e. in }~~B.
	\end{array}
	\end{equation}
	We claim that $w_{0}\ne 0$. Suppose, by contradiction, $w_{0}= 0$. From the definition
	of $\mathcal{N}$ and \eqref{eq:5.14} , we have that $ \displaystyle\lim_{n \rightarrow +\infty} \Vert w_n \Vert^2 =0$, which contradicts Lemma \ref{lem3}. Hence, $w_{0}\ne 0$ .
	
	From the lower semi continuity of norm and \eqref{eq:5.14}, it follows that
 	\begin{equation}\label{eq:5.15}
 	\Vert w_{0}\Vert^2 \leq  \liminf_{n \rightarrow +\infty}(\Vert w_n\Vert^2 )
 	\end{equation}	

 On the other hand, by using $\langle J'(w_n),w_n\rangle =0$ and (\ref{eq:5.14}), we have

 	\begin{equation}\label{eq:5.16} \liminf_{n \rightarrow +\infty}(\Vert w_n\Vert^2 ) =\liminf_{n \rightarrow +\infty} \int_{B}f(x, w_n) w_n dx=\int_{B}f(x, w_{0}) w_{0} dx.\end{equation}
 From \eqref{eq:5.15} and \eqref{eq:5.16} we deduce that $\langle \mathcal{J}'(w_{0}),w_{0}\rangle \leq 0$.
 Then, as in Lemma \ref{lem2.2} this implies that there exists $s\in (0, 1] $ such that $sw_{0}\in\mathcal{N}$. Thus, by the lower semi continuity of norm, Lemma \ref{lem11} and Lemma \ref{lemma9}, we get that

$$\begin{array}{rclll}
\displaystyle m \leq \mathcal{J}(sw_{0})&=&\mathcal{J}(sw_{0}) -\displaystyle\frac{1}{2} \langle \mathcal{J}'(sw_{0}),sw_{0}\rangle\\&=&\displaystyle\frac{1}{2}\int_{B}\big(f(x,sw_{0})sw_{0}-2F(x,sw_{0})\big)dx
\\ &\leq& \displaystyle\frac{1}{2}\int_{B}\big(f(x,w_{0})w_{0}-2F(x,w_{0})\big)dx \\
&<&\displaystyle\frac{1}{2}\int_{B}\big(f(x,w_{0})w_{0}-2F(x,w_{0})\big)dx \\
&\leq&\displaystyle\liminf_{n\rightarrow +\infty}\Big[\displaystyle\frac{1}{2}\|w_{n}\|^{2} -\int_{B}F(x,w_{n})dx-\frac{1}{2}\|w_{n}\|^{2}+\frac{1}{2}\int_{B}f(x,w_{n})w_{n}dx\Big]\\
&\leq&\displaystyle \liminf_{n\rightarrow+\infty}\big[\mathcal{J}(w_n) -\displaystyle\frac{1}{2} \langle \mathcal{J}'(w_n), w_n\rangle\big] = m.
\end{array}$$	
Therefore, we get that $\mathcal{J}(sw_{0})=m$, which is the desired conclusion.

\textbf{Proof of Theorem \ref{th1.3}. }From Lemma \ref{lem12} there exists $w_{0}$ such that $\mathcal{J}(w_{0})=m$. Now, by Lemma \ref{lem5}, we deduce that $\mathcal{J}'(w_{0})=0$.
So, $w_{0}$ is a solution for our problem $(P)$ .

\end{document}